\documentclass[12pt, reqno]{amsart}
\allowdisplaybreaks
\begin{document}
\setcounter{page}{1}

\title[\hfilneg \hfil A generalized Gronwall ]
{A generalized Gronwall Inequality for Caputo Fractional Dynamic delta operator }

\author[Deepak Pachpatte \hfil \hfilneg]
{Deepak B. Pachpatte}

\address{Deepak B. Pachpatte \newline
 Dept. of Mathematics,
 Dr. B. A. M. University, Aurangabad,
 Maharashtra 431004, India}
\email{pachpatte@gmail.com}


\begin{abstract}
  In this paper we obtain generalized Gronwall type inequality using Caputo Fractional delta operator. Also we have  obtained the existence of solution of Cauchy's Type problem on fractional dynamic equations using dynamic delta operator. Applying the obtained  inequality we study the properties of solution on fractional dynamic equations.

\end{abstract}

\maketitle

\section{Introduction}

	Fractional Calculus is an important tool which generalizes the differential and integral calculus of arbitrary order. In this it is possible to define the differentiation and integration for non-integer order. Fractional calculus is more suitable for modeling the real world problems in various branches of science and engineering. In year 1989 Stefan Hilger introduced time scale calculus  which is  unification of the differential and difference calculus \cite{HIG}. Since then many authors have studied in properties and various applications of dynamic equations on time scales \cite{Geo}.

	 In \cite{Adj, ALA, Ye} authors studied the  Gronwall type inequality and its applications on fractional Differential Equations using various fractional operators.
	 On the other hand \cite{AHM, Ben, Yan, Zhu1,Zhu2, Zha1} authors have combined the fractional calculus and time scale calculus and obtained results on existence and some properties of fractional differential equations on time scales. Basic information on time scale calculus can be found in \cite{Boh1, Boh2}

	The  basic theory  on fractional dynamic calculus and equations on time scales can be found \cite{Ana1, Ana2, Geo1}. This types of problems has applications in studying the properties of various processes in materials \cite{Bal}.
	Inspired  by above literature in this paper we obtain the estimates on Gronwall type inequality and obtain the existence of solution Cauchy's Type problem on fractional dynamic equations on time scales. Using the  obtained  inequality we study the properties of Cauchy's type of problem such as continuous dependence of solution.
	
\section{Preliminaries}
Now in this section we give some basic definitions and theorems  which are used in our subsequent discussions.

We denote by $\mathbb{T}$ any time scale which has a topology which it inherits from standard topology on $\mathbb{T}$. We denote $C_{rd}$ for the set of all rd-continuous functions. For more basic information  on time scale see \cite{Boh1, Boh2}.

Now as in  \cite{Kul} we construct the metric space where $\left[ {t_0 ,\infty } \right)_\mathbb{T}  = I_\mathbb{T} $.Now consider the space  function $C_{rd} \left( {I_\mathbb{T} ,\mathbb{R}^n } \right)$ such that $\mathop {\sup }\limits_{t \in I_\mathbb{T} } \frac{{u\left( t \right)}}{{e_\eta  \left( {t,t_0 } \right)}} < \infty$ where $\eta  > 0$. This space is denoted by $C_{rd}^\eta  \left( {I_\mathbb{T} ,\mathbb{R}^n } \right)$.

We couple the space $C_{rd}^\eta  \left( {I_\mathbb{T} ,\mathbb{R}^n } \right)$ by suitable metric
\[
m_\eta ^\infty  \left( {u,v} \right) = \mathop {\sup }\limits_{t \in I_T } \frac{{\left| {u\left( t \right) - v\left( t \right)} \right|}}{{e_\eta  \left( {t,t_0 } \right)}},
\]
where the norm is defined as
\[
\left| u \right|_\eta ^\infty   = \mathop {\sup }\limits_{t \in I_\mathbb{T} } \frac{{\left| {u\left( t \right)} \right|}}{{e_\eta  \left( {t,t_0 } \right)}}.
\]
More properties of $m_\eta ^\infty$ and $\left| . \right|_\eta ^\infty $  can be found in \cite{Kul}.

  We define delta power function as
\paragraph{\textbf{Definition 2.1} \cite{Geo1}}. Let $\alpha \in \mathbb{R}$, we define the generalized delta power function $h_{\alpha}$ on $\mathbb{T}$
as follows:
\[
h_\alpha  \left( {t,t_0 } \right) = L^{ - 1} \left( {\frac{1}{{z^{\alpha  + 1} }}} \right)\left( t \right),\,\,\,\,\,t \ge t_0 ,
\]
for all $z \in C\backslash \{ 0\} $ such that $L^{ - 1}$ exists, $t \ge t_0$. The fractional generalized delta power function $h_\alpha(t,s)$ on $T$, $t \ge s \ge t_0 $  which is defined as the shift of $h_\alpha  \left( {t,t_0 } \right)$ i.e.,
\[
h_\alpha  \left( {t,s} \right) = \,\,\,\,\widehat{h_\alpha ({.,t_0 })}\,( {t,s} ),\,\,\,\,\,\,\,t,s \in T,\,\,\,\,\,t \ge s \ge t_0 \,.
\]
Now we define the Riemann-Liouville Fractional delta integral and Riemann Liouville Fractional delta derivative as follows. Suppose $\alpha \ge 0$ and $[-\overline \alpha]$ denote the integral part of $- \alpha$.
\paragraph{\textbf{Definition 2.2} \cite{Geo1,Zhu1}}
For a function $f:\mathbb{T}  \rightarrow \mathbb{R}$ the Riemann Liouville fractional delta integral of order $\alpha$ defined by
\[
I_{\Delta ,t_0 }^0 f\left( t \right) = f(t),
\]
\begin{align*}
\left( {I_{\Delta ,t_0 }^0 f} \right)\left( t \right)
& = \left( {h_{\alpha  - 1} \left( {.,t_0 } \right)*f} \right)(t) \\
&= \int\limits_{t_0 }^t \widehat{h_{\alpha  - 1} \left( {.,t_0 } \right)} \left( {t,\sigma \left( u \right)} \right)f\left( u \right)\Delta u \\
&= \int\limits_{t_0 }^t {h_{\alpha  - 1} } \left( {t,\sigma \left( u \right)} \right)f\left( u \right)\Delta u.
\end{align*}

\paragraph{\textbf{Definition 2.3} \cite{Geo1,Zhu1}}
Let $\alpha \ge 0$, $m=-\overline[-\alpha]$, $f:\mathbb{T}  \rightarrow \mathbb{R}$. For $s,t \in \mathbb{T}^{k^m }$, $s<t$, the Riemann-Liouville fractional delta derivative of order $\alpha$ is defined by the expression
\[
D_{\Delta ,s}^\alpha  f\left( t \right) = D_\Delta ^m I_{\Delta ,s}^{m - \alpha } f\left( t \right),\,\,\,\,\,\,t \in \mathbb{T},
\]
if it exists. For $\alpha <0$ we define
\[
D_{\Delta ,s}^\alpha  f\left( t \right) = I_{\Delta ,s}^{ - \alpha } f\left( t \right),\,\,\,\,\,t,s \in T,\,\,\,t > s,
\]
\[
I_{\Delta ,s}^\alpha  f\left( t \right) = D_{\Delta ,s}^{ - \alpha } f\left( t \right),\,\,\,\,t,s \in T^{k^m } ,\,\,\,t > s,\,\,\,r = \overline {\left[ { - \alpha } \right]}  + 1.
\]

\paragraph{}Now we define Caputo Fractional delta derivative as
\paragraph{\textbf{Definition 2.4} \cite{Geo1}}
Let $t \in T$. The Caputo fractional delta derivative of order $\alpha \ge 0$ is defined via the Riemann-Liouville fractional delta derivative as follows:
\[
{}^CD_{\Delta ,t_0 }^\alpha  f\left( t \right) = D_{\Delta ,t_0 }^\alpha  \left( {f\left( t \right) - \sum\limits_{k = 0}^{m - 1} {h_k \left( {t,t_0 } \right)} f^{\Delta ^k } \left( {t_0 } \right)} \right),\,\,\,\,\,\,\,\,t > t_0 ,
\]
where $m = \overline {\left[ \alpha  \right]}  + 1$ if $\alpha  \notin \mathbb{N}$, $m = \overline {\left[ \alpha  \right]}$ if $\alpha  \in \mathbb{N}$.

\section{Gronwall Type Inequality}
Now we give the Gronwall type inequality using Caputo Fractional delta operator and we prove this by iteration. Suppose $\alpha \ge 0$ and $[-\overline \alpha]$ denote the integral part of $- \alpha$.

\paragraph{\textbf{Theorem 3.1}}
Let $\alpha>0$, $y,u:\mathbb{T}\rightarrow \mathbb{R}$ be two non-negative integrable functions and $v$ be non negative, non decreasing and rd-continuous function,  $v(t) \le B$ be a constant. If
\[
y\left( t \right) \le u\left( t \right) + v\left( t \right)\int\limits_{t_0 }^t {h_{\alpha  - 1} \left( {t,\sigma \left( \tau  \right)} \right)} \Delta \tau ,
\tag{3.1}\]
then
\[
y\left( t \right) \le u\left( t \right) + \int\limits_{t_0 }^t {\left[ {\sum\limits_{k = 1}^\infty  {\left( {v\left( t \right)} \right)^k h_{k\alpha  - 1} \left( {t,\sigma \left( \tau  \right)} \right)u\left( \tau  \right)} } \right]} \Delta \tau.
\tag{3.2}\]
\paragraph{\textbf{Proof}}
 Define a function $Q$ by
\[
Q\psi \left( t \right) = v\left( t \right)\int\limits_{t_0 }^t {h_{\alpha  - 1} } \left( {t,\sigma \left( \tau  \right)} \right)\psi \left( \tau  \right)\Delta \tau,
\tag{3.3}\]
then we get
\[
y\left( t \right) \le u\left( t \right) + Qy\left( t \right).
\tag{3.4}\]
Now taking iteration of $(3.4)$ consecutively we get for $n \in N$
\[
y(t) \le \sum\limits_{k = 0}^{n - 1} {Q^k u\left( t \right)}  + Q^n u\left( t \right).
\tag{3.5}\]
Now we prove by induction hypotheses that if $\psi$ is non negative function then
\[
Q^k \psi \left( t \right) \le \int\limits_{t_0 }^t {\left( {v\left( t \right)} \right)^k h_{k\alpha  - 1} } \left( {t,\sigma \left( s \right)} \right)\psi \left( s \right)\Delta s.
\tag{3.6}\]
If $k=1$ the result is obvious. Suppose the formula is valid for $k \in N$. We have
\begin{align*}
&Q^{k + 1} \psi \left( t \right)\\
&= Q.Q^k \psi \left( t \right) \\
& \le v\left( t \right)\int\limits_{t_0 }^t {h_{\alpha  - 1} \left( {t,\sigma \left( \tau  \right)} \right)} \left[ {\int\limits_{t_0 }^\tau  {\left( {v\left( \tau  \right)} \right)^k h_{k\alpha  - 1} } \left( {\tau ,\sigma \left( s \right)} \right)\psi \left( s \right)\Delta s} \right]\Delta \tau.  \\
\tag{3.7}
\end{align*}
We have $v$  non decreasing $v(\tau) \le v(t)$ for $\tau \le t$, from $(3.7)$
\begin{align*}
&Q^{k + 1} \psi \left( t \right) \\
&\le \left( {v\left( t \right)} \right)^{k + 1} \int\limits_{t_0 }^t {\left[ {\int\limits_{t_0 }^\tau  {h_{\alpha  - 1} \left( {t,\sigma \left( \tau  \right)} \right)h_{k\alpha  - 1} \left( {\tau ,\sigma \left( s \right)} \right)\Delta \tau } } \right]} \psi \left( s \right)\Delta s.
\tag{3.8}
\end{align*}
From \cite{Ana2} and properties of the inner integral is
\[
\int\limits_{t_0 }^t {h_{\alpha  - 1} \left( {t,\sigma \left( \tau  \right)} \right)h_{k\alpha  - 1} \left( {\tau ,\sigma \left( s \right)} \right)\Delta \tau }  = h_{(k + 1)\alpha  - 1} \left( {\tau ,\sigma \left( s \right)} \right).
\tag{3.9}\]
Then from $(3.8)$ we get
\[
Q^{k + 1} \psi \left( t \right) \le \left( {v\left( t \right)} \right)^{k + 1} \int\limits_{t_0 }^t {h_{(k + 1)\alpha  - 1} \left( {\tau ,\sigma \left( s \right)} \right)} \psi \left( s \right)\Delta s.
\tag{3.10}\]
This proves that
\[
Q^n \psi \left( t \right) \le \int\limits_{t_0 }^t {\left( {v\left( t \right)} \right)^k h_{k\alpha  - 1} \left( {\tau ,\sigma \left( s \right)} \right)} \psi \left( s \right)\Delta s.
\tag{3.11}\]
Now we show that $\psi ^n y\left( t \right) \to 0$ as $ n \to \infty$.
Since  $g(t)$ is rd-continuous and there exists $B>0$ such that $g(t) \le B$ then we have
\[
Q^n y\left( t \right) \le \int\limits_{t_0 }^t {B^N h_{N\alpha  - 1} \left( {\tau ,\sigma \left( s \right)} \right)} y\left( s \right)\Delta s,
\tag{3.12}\]
$\rightarrow 0$ as $n \to \infty $.

Therefore we have from $(3.5)$
\[
y\left( t \right) \le \sum\limits_{k = 0}^\infty  {Q^k f\left( t \right)} .
\tag{3.13}\]
Thus we get
\begin{align*}
 y\left( t \right)
& \le \sum\limits_{k = 0}^\infty  {Q^k f\left( t \right)}  \\
& \le u(t) + \int\limits_{t_0 }^t {\sum\limits_{k = 1}^\infty  {\left( {v\left( t \right)} \right)^k h_{k\alpha  - 1} \left( {\tau ,\sigma \left( s \right)} \right)f\left( t  \right)} } \Delta \tau,
\tag{3.14}
\end{align*}
for $t \in \mathbb{T}$, which is required inequality.

\section{Existence and Uniqueness}

Now we consider the Cauchy's type of problem with Caputo fractional delta derivative, suppose $\alpha >0$
\[
{}^CD_{\Delta ,t_0 }^\alpha  u\left( t \right) = f\left( {t,u\left( t \right)} \right),\,\,\,\,\,\,t \in I_\mathbb{T} ,
\tag{4.1}\]

with the initial condition
\[
{}^CD_{\Delta ,t_0 }^\alpha  u\left( {t_0 } \right) = \overline w ,
\tag{4.2}\]
where $f:\mathbb{T}\times \mathbb{R} \rightarrow \mathbb{R}$ is a given function and $0 <\alpha <1$.

Let $L_{\Delta}[t_0,a)$ denote the space of $\Delta$ Lebesgue summable function in $[t_0,a)$. Define the space
\[
L_\Delta ^\alpha  \left[ {t_0 ,\alpha } \right) = \left\{ {y \in L_\Delta  \left[ {t_0 ,a} \right):D_{\Delta ,t_0 }^\alpha  y \in L_\Delta  \left[ {t_0 ,a} \right)} \right\}.
\]
Then from Theorem $52$, \cite{Zhu1}, $(4.1)$ and $(4.2)$ is equivalent to
\[
u\left( t \right) = w h_{\alpha  - 1} \left( {t,t_0 } \right) + \int\limits_{t_0 }^t {h_{\alpha  - 1} \left( {t,\sigma \left( \tau  \right)} \right)f\left( {\tau ,u\left( \tau  \right)} \right)\Delta \tau }.
\tag{4.3}\]
 Now we give the existence of solution in next theorem
\paragraph{\textbf{Theorem 4.1}}. Let $L \ge 0$  be a constant. Suppose the function $f$ in $(4.1)$ be rd-continuous and satisfy
\[
\left| {f\left( {x_1 ,x_2 } \right) - f\left( {x_1 ,\overline {x_2 } } \right)} \right| \le L\left| {x_2  - \overline {x_2 } } \right|,
\tag{4.4}\]
and let
\[
p_1  = \mathop {\sup }\limits_{t \in I_T } \frac{1}{{e_\eta  \left( {t,t_0 } \right)}}\left| {wh_{\alpha  - 1} \left( {t,t_0 } \right) + \int\limits_{t_0 }^t {h_{\alpha  - 1} \left( {t,t_0 } \right)f\left( {\tau ,0} \right)\Delta \tau } } \right| < \infty .
\tag{4.5}\]
If $\frac{L}{\eta } < 1$ then equation $(4.1)$ has a unique solution $u \in C_{rd}^\eta  \left( {I_\mathbb{T} ,\mathbb{R}^n } \right) $.
\paragraph{\textbf{Proof.}} Let $u \in C_{rd}^\eta  \left( {I_\mathbb{T} ,\mathbb{R}^n } \right)$ and define operator $G$ by
\[
\left( {Gu} \right)\left( t \right) = w\,h_{\alpha  - 1} \left( {t,t_0 } \right) + \int\limits_{t_0 }^t {h_{\alpha  - 1} } \left( {t,\sigma \left( \tau  \right)} \right)f\left( {\tau ,u\left( \tau  \right)} \right)\Delta \tau ,
\tag{4.6}\]
for $t \in I_\mathbb{T} $.

We prove that $G$ maps $C_{rd}^\eta  \left( {I_\mathbb{T} ,\mathbb{R}^n } \right)$ into itself and is a contraction map.

From (4.6) we have
\begin{align*}
\left( {Gu} \right)\left( t \right)
& = w\,h_{\alpha  - 1} \left( {t,t_0 } \right) + \int\limits_{t_0 }^t {h_{\alpha  - 1} } \left( {t,\sigma \left( \tau  \right)} \right)f\left( {\tau ,u\left( \tau  \right)} \right)\Delta \tau  \\
& - \int\limits_{t_0 }^t {h_{\alpha  - 1} } \left( {t,\sigma \left( \tau  \right)} \right)f\left( {\tau ,0} \right)\Delta \tau  + \int\limits_{t_0 }^t {h_{\alpha  - 1} } \left( {t,\sigma \left( \tau  \right)} \right)f\left( {\tau ,0} \right)\Delta \tau.
\tag{4.7}
\end{align*}
We prove that $G$ maps $C_{rd}^\eta  \left( {I_\mathbb{T} ,\mathbb{R}^n } \right)$ into itself and is a contraction map. From $(4.7)$ we have
\begin{align*}
\left| {Gu} \right|_\eta ^\infty
& = \mathop {\sup }\limits_{t \in I_\mathbb{T} } \frac{{\left| {\left( {Gu} \right)\left( t \right)} \right|}}{{e{}_\eta \left( {t,t_0 } \right)}} \\
& = \mathop {\sup }\limits_{t \in I_\mathbb{T} } \frac{1}{{e{}_\eta \left( {t,t_0 } \right)}}\left| {w\,h_{\alpha  - 1} \left( {t,t_0 } \right) + \int\limits_{t_0 }^t {h_{\alpha  - 1} } \left( {t,\sigma \left( \tau  \right)} \right)f\left( {\tau ,u\left( \tau  \right)} \right)\Delta \tau } \right. \\
&\left. { - \int\limits_{t_0 }^t {h_{\alpha  - 1} } \left( {t,\sigma \left( \tau  \right)} \right)f\left( {\tau ,0} \right)\Delta \tau  + \int\limits_{t_0 }^t {h_{\alpha  - 1} } \left( {t,\sigma \left( \tau  \right)} \right)f\left( {\tau ,0} \right)\Delta \tau } \right| \\
& \le \mathop {\sup }\limits_{t \in I_\mathbb{T} } \frac{1}{{e{}_\eta \left( {t,t_0 } \right)}}\left| {w\,h_{\alpha  - 1} \left( {t,t_0 } \right) + \int\limits_{t_0 }^t {h_{\alpha  - 1} } \left( {t,\sigma \left( \tau  \right)} \right)f\left( {\tau ,0} \right)\Delta \tau } \right| \\
& + \mathop {\sup }\limits_{t \in I_\mathbb{T} } \frac{1}{{e{}_\eta \left( {t,t_0 } \right)}}\left| {\int\limits_{t_0 }^t {h_{\alpha  - 1} } \left( {t,\sigma \left( \tau  \right)} \right)f\left( {\tau ,u\left( \tau  \right)} \right)\Delta \tau } \right. \\
&\left. { - \int\limits_{t_0 }^t {h_{\alpha  - 1} } \left( {t,\sigma \left( \tau  \right)} \right)f\left( {\tau ,0} \right)\Delta \tau } \right| \\
&= p_1  + \mathop {\sup }\limits_{t \in I_\mathbb{T} } \frac{1}{{e{}_\eta \left( {t,t_0 } \right)}}\int\limits_{t_0 }^t {h_{\alpha  - 1} } \left( {t,\sigma \left( \tau  \right)} \right)\left| {f\left( {\tau ,u\left( \tau  \right)} \right) - f\left( {\tau ,0} \right)} \right|\Delta \tau \\
&= p_1  + \mathop {\sup }\limits_{t \in I_\mathbb{T} } \frac{1}{{e{}_\eta \left( {t,t_0 } \right)}}\int\limits_{t_0 }^t {h_{\alpha  - 1} } \left( {t,\sigma \left( \tau  \right)} \right)L(u(\tau ))\Delta \tau  \\
&= p_1  + L\left| u \right|_\eta ^\infty  \mathop {\sup }\limits_{t \in I_\mathbb{T} } \frac{1}{{e{}_\eta \left( {t,t_0 } \right)}}\int\limits_{t_0 }^t {h_{\alpha  - 1} } \left( {t,\sigma \left( \tau  \right)} \right)e{}_\eta \left( {\tau ,t_0 } \right)\Delta \tau  \\
& = p_1  + L\left| u \right|_\eta ^\infty  \mathop {\sup }\limits_{t \in I_\mathbb{T} } \frac{1}{{e{}_\eta \left( {t,t_0 } \right)}}I_{t_0 }^\Delta  \left( {e{}_\eta \left( {\tau ,t_0 } \right)} \right)\Delta \tau  \\
&\le p_1  + L\left| u \right|_\eta ^\infty  \mathop {\sup }\limits_{t \in I_\mathbb{T} } \frac{1}{{e{}_\eta \left( {t,t_0 } \right)}}\left( {\frac{{e{}_\eta \left( {t,t_0 } \right) - 1}}{\eta }} \right) \\
&\le p_1  + L\left| u \right|_\eta ^\infty  \frac{1}{\eta }\left( {1 - \frac{1}{{e{}_\eta \left( {t,t_0 } \right)}}} \right) \\
&= p_1  + \frac{L}{\eta }\left| u \right|_\eta ^\infty   \\
&< \infty.
\tag{4.8}
\end{align*}
This proves that $G$ maps $C_{rd}^\eta  \left( {I_\mathbb{T} ,\mathbb{R}^n } \right)$ into itself.

Now we prove that $G$ is a Contraction map.

Let $x,y \in C_{rd}^\eta  \left( {I_\mathbb{T} ,\mathbb{R}^n } \right)$ then from $(3.7)$ and by hypotheses we get
\begin{align*}
m_\eta ^\infty  \left( {Gx,Gy} \right)
&= \mathop {\sup }\limits_{t \in I_\mathbb{T} } \frac{{\left| {\left( {Gx} \right)\left( t \right) - \left( {Gy} \right)\left( t \right)} \right|}}{{e{}_\eta \left( {t,t_0 } \right)}} \\
&= \mathop {\sup }\limits_{t \in I_\mathbb{T} } \frac{1}{{e{}_\eta \left( {t,t_0 } \right)}}\left| {\int\limits_{t_0 }^t {h_{\alpha  - 1} } \left( {t,\sigma \left( \tau  \right)} \right)f\left( {\tau ,x\left( \tau  \right)} \right)\Delta \tau } \right. \\
&\left. { - \int\limits_{t_0 }^t {h_{\alpha  - 1} } \left( {t,\sigma \left( \tau  \right)} \right)f\left( {\tau ,y\left( \tau  \right)} \right)\Delta \tau } \right| \\
&\le \mathop {\sup }\limits_{t \in I_\mathbb{T} } \frac{1}{{e{}_\eta \left( {t,t_0 } \right)}}\left| {\int\limits_{t_0 }^t {h_{\alpha  - 1} } \left( {t,\sigma \left( \tau  \right)} \right)L\frac{{\left| {x\left( \tau  \right) - y\left( \tau  \right)} \right|}}{{e{}_\eta \left( {t,t_0 } \right)}}e{}_\eta \left( {t,t_0 } \right)} \right|\Delta \tau \\
&= \mathop {\sup }\limits_{t \in I_\mathbb{T} } \frac{1}{{e{}_\eta \left( {t,t_0 } \right)}}\left| {\int\limits_{t_0 }^t {h_{\alpha  - 1} } \left( {t,\sigma \left( \tau  \right)} \right)m_\eta ^\infty  \left( {x,y} \right)e{}_\eta \left( {t,t_0 } \right)\Delta \tau } \right| \\
& = \mathop {\sup }\limits_{t \in I_\mathbb{T} } \frac{1}{{e{}_\eta \left( {t,t_0 } \right)}}Lm_\eta ^\infty  \left( {x,y} \right)\int\limits_{t_0 }^t {h_{\alpha  - 1} } \left( {t,\sigma \left( \tau  \right)} \right)e{}_\eta \left( {t,t_0 } \right)\Delta \tau  \\
& = Lm_\eta ^\infty  \left( {x,y} \right)\mathop {\sup }\limits_{t \in I_\mathbb{T} } \frac{1}{{e{}_\eta \left( {t,t_0 } \right)}}\left( {\frac{{e{}_\eta \left( {t,t_0 } \right) - 1}}{\eta }} \right) \\
&\le \frac{L}{\eta }m_\eta ^\infty  \left( {x,y} \right).
\tag{4.9}
\end{align*}
Since $\frac{L}{\eta }<1$. Thus $G$ has a unique fixed point in $C_{rd}^\eta  \left( {I_\mathbb{T} ,\mathbb{R}^n } \right)$ from Banach Fixed point theorem. The fixed point of G is a solution of equation $(4.1)$. This completes the proof of theorem.

\section{Continuous Dependence}
In this section we obtain the results for continuous dependence of solution of $(4.1)$.
Now consider the equation $(4.1)$ and the corresponding equation
\[
{}^CD_{\Delta ,t_0 }^\alpha  v\left( t \right) = \overline f \left( {t,v\left( t \right)} \right)\,,\,t \in I_\mathbb{T}, \,\,
\tag{5.1}\]
with initial condition
\[
{}^CD_{\Delta ,t_0 }^\alpha  v\left( t \right) = \overline w,
\tag{5.2}\]
where $f:I_\mathbb{T} \, \times \mathbb{R} \to \mathbb{R}$ and $\overline w $ is a given constant.
\paragraph{} Now we give the theorem which deals with continuous dependence of solution of $(4.1)$.
\paragraph{\textbf{Theorem 5.1}}
Suppose the function $f$ in $(4.1)$ rd-continuous and satisfy the condition $(4.4)$. Let $v(t)$ be solution of $(5.1)$ and

\begin{align*}
 H(t)
&= \left| {wh_{\alpha  - 1} \left( {t,t_0 } \right) - \overline w h_{\alpha  - 1} \left( {t,t_0 } \right)} \right| \\
&+ \left| {\int\limits_{t_0 }^t {h_{\alpha  - 1} \left( {t,\sigma \left( \tau  \right)} \right)f\left( {\tau ,v\left( \tau  \right)} \right)\Delta \tau } } \right. \\
&\left. { - \int\limits_{t_0 }^t {h_{\alpha  - 1} \left( {t,\sigma \left( \tau  \right)} \right)\overline f \left( {\tau ,v\left( \tau  \right)} \right)\Delta \tau } } \right|,
\tag{5.3}
\end{align*}
where $f$ and $\overline f$ are functions in $(4.1)$ and $(5.1)$. Then the solution $u(t), t \in I_\mathbb{T}$ of $(4.1)$ dependence on functions on right hand side of $(4.1)$ and
\[
\left| {u(t) - v(t)} \right| \le H\left( t \right) + \int\limits_{t_0 }^t {\left[ {\sum\limits_{k = 1}^\infty  {L^k h_{k\alpha  - 1} \left( {t,\sigma \left( \tau  \right)} \right)H\left( \tau  \right)} } \right]} \Delta \tau.
\tag{5.4}\]
\paragraph{\textbf{Proof.}} The solutions of the equation $(4.1)-(4.2)$ and $(5.1)-(5.2)$ are
\[
u\left( t \right) = w\,h_{\alpha  - 1} \left( {t,t_0 } \right) + \int\limits_{t_0 }^t {h_{\alpha  - 1} \left( {t,\sigma \left( \tau  \right)} \right)f\left( {\tau ,u\left( \tau  \right)} \right)} \Delta \tau ,
\tag{5.5}\]
and
\[
v\left( t \right) = \overline w \,h_{\alpha  - 1} \left( {t,t_0 } \right) + \int\limits_{t_0 }^t {h_{\alpha  - 1} \left( {t,\sigma \left( \tau  \right)} \right)f\left( {\tau ,v\left( \tau  \right)} \right)} \Delta \tau ,
\tag{5.6}\]
respectively.

We have
\begin{align*}
\left| {u\left( t \right) - v\left( t \right)} \right|
&\le \left| {w\,h_{\alpha  - 1} \left( {t,t_0 } \right) - \overline w \,h_{\alpha  - 1} \left( {t,t_0 } \right)} \right| \\
& + \left| {\int\limits_{t_0 }^t {h_{\alpha  - 1} \left( {t,\sigma \left( \tau  \right)} \right)f\left( {\tau ,u\left( \tau  \right)} \right)} \Delta \tau } \right. \\
&\left. { - \int\limits_{t_0 }^t {h_{\alpha  - 1} \left( {t,\sigma \left( \tau  \right)} \right)f\left( {\tau ,v\left( \tau  \right)} \right)} \Delta \tau } \right| \\
 & + \left| {\int\limits_
{t_0 }^t {h_{\alpha  - 1} \left( {t,\sigma \left( \tau  \right)} \right)f\left( {\tau ,v\left( \tau  \right)} \right)} \Delta \tau } \right. \\
&\left. { - \int\limits_{t_0 }^t {h_{\alpha  - 1} \left( {t,\sigma \left( \tau  \right)} \right)\overline f \left( {\tau ,v\left( \tau  \right)} \right)} \Delta \tau } \right| \\
& \le H(t) + \int\limits_{t_0 }^t {h_{\alpha  - 1} \left( {t,\sigma \left( \tau  \right)} \right)L\left| {u\left( \tau  \right) - v\left( \tau  \right)} \right|} \Delta \tau .
\tag{5.7}
\end{align*}
Now an application of Theorem $(3.1)$ to equation $(5.7)$ yields the required inequality $(5.6)$.


\begin{thebibliography}{999}

\bibitem {Adj} Y. Adjabi, F. Jarad and T Abdeljawad,
\newblock  On Generalized Fractional Operators and a Gronwall Type inequality With Applications,
\newblock {\em Filomat} 31:17(2017), 5457-5473.

\bibitem {AHM} A. Ahamadkhanlu and M. Jahanshahi,
\newblock On the existence and uniqueness of solution of initial value problem for fractional order differential equations on time scales,
\newblock {\em Bull. Iranian. Math. Soc.}, Vol 38(1),2012, pp. 241-252.

\bibitem {Ana1} G. Anastassiou,
\newblock Frontiers in time scales and Inequalities,
\newblock {\em World Scientific Publishing Company}, (2015).

\bibitem {Ana2} G. Anastassiou,
\newblock Principle of delta fractional calculus on time scales and inequalities,
\newblock {\em Math. Comput. Modelling}, 52(2010), 556-566.

\bibitem {ALA} R. Almedia,
\newblock A Gronwall inequality for a General Caputo Fractional Operator,
\newblock {\em Math. Inequal. Appl.}, 20(4),(2017), pp 1089-1105.

\bibitem {Bal} D. Balenu, K. Diethelm, E. Scales and J. Trujillo,
\newblock Fractional Calculus Models and Numerical Methods,
\newblock {\em World Scientific Publishing Company}, (2017).


\bibitem {Ben} N. Benkhettou, A. Hammoudi and  D.Torres,
\newblock Existence and uniqueness of solution for a fractional Riemann, liouville initial value problem on time scales,
\newblock {\em J King Saud Univ Sci},Vol 28, Issue 1, January 2016, pp 87-92

\bibitem {Boh1} M. Bohner  and A. Peterson,
\newblock Dynamic equations on time scales,
\newblock {\em Birkhauser Boston/Berlin}, (2001).

\bibitem {Boh2} M. Bohner  and A. Peterson,
\newblock Advances in Dynamic equations on time scales,
\newblock {\em Birkhauser Boston/Berlin}, (2003).

\bibitem {Geo} Sevtlin G. Georgiev,
\newblock Integral equations on time scales,
\newblock {\em Atlantis Press}, 2016

\bibitem {Geo1} Sevtlin G. Georgiev,
\newblock Fractional dynamic calculus and fractional dynamic equations on time scales,
\newblock {\em Springer}, 2017

\bibitem {HIG} S. Hilger,
\newblock Analysis on Measure chain-A unified approach to continuous discrete calculus,
\newblock {\em Results. Math.}, 18:18-56, 1990.


\bibitem {Kul} T. Kulik and C. C. Tisdell,
\newblock Volterra integral equations on time scales: Basic qualitative and quantitative results with applications to initial value problems on unbounded domains,
\newblock {\em Int. J. Difference Equ.}, Vol. 3, No 1 (2008), 103-133.

\bibitem {Yan} R.A. Yan, S.R. Sun and Z. L. Han
\newblock Existence of solutions of boundary value problems for Caputo Fractional Differential equations on time scales,
\newblock {\em Bull. Iranian. Math. Soc.}, Vol 42(2016), No.2, pp. 247-262

\bibitem {Ye} H. Ye, J. Gao and Y. Ding,
\newblock A generalized Gronwall inequality and its application to a fractional differential equation,
\newblock {\em J. Math. Anal. Appl.}, 328(2007), 1075-1081.

\bibitem {Zhu1} Jiang Zhu and Ying Zhu,
\newblock Fractional Cauchy Problem with Riemann-Liouville Fractional Delta derivative on Time Scales,
\newblock {\em Abstr. Appl. Anal.}, Art. Id 401596, Vol 2013 , pp1-19.

\bibitem {Zhu2} Jiang Zhu and Ling Wu,
\newblock Fractional Cauchy Problem with Caputo Nabla derivative on Time Scales,
\newblock {\em Abstr. Appl. Anal.}, Art. Id 486054, Vol 2015 , pp1-23.

\bibitem {Zha1} Xiaozhi Zhang and Chuanxi Zhu,
\newblock Cauchy problem for a class of fractional differential equations on time scales,
\newblock {\em  Int. J. Comput. Math.}, Vol 91, No 3, 527-538.


 \end{thebibliography}
\end{document}